\documentstyle[12pt]{article}
\catcode`\@=11
\@addtoreset{equation}{section}

\catcode`\@=12
\newtheorem{Theorem}{Theorem}[section]
\newtheorem{Definition}{Definition}[section]

\newtheorem{Lemma}{Lemma}[section]
\newtheorem{Corollary}{Corollary}[section]

\title{Proofs On Arnold Chord Conjecture and Weinstein Conjecture in 
$M\times C$\thanks{Project 19871044 Supported by NSF}}

\author{Renyi Ma\\
Department of Mathmatics \\
Tsinghua University \\
Beijing, 100084\\
People's Republic of China}

\date { }

\begin{document}
\textwidth=125mm
\textheight=185mm
\parindent=8mm
\frenchspacing
\maketitle

\begin{abstract}
In this article, we give new proofs on  
the some cases on Arnold chord conjecture and Weinstein conjecture in $M\times C$ which 
includes the previous works as special cases. 
\end{abstract}


\noindent{\bf 2000MR Subject Classification} 32Q65,53D35,53D12

\section{Introduction and results}

\subsection{Arnold chord conjecture}

Let $\Sigma$ be a smooth closed oriented manifold of dimension
$2n-1$. A contact form on $\Sigma$ is a $1-$form such that
$\lambda \wedge (d\lambda )^{n-1}$ is a volume form on $\Sigma$.
Associated to $\lambda$ there is the so-called Reed vectorfield $X_\lambda $ defined
by $i_X\lambda   \equiv 1, \ \ i_Xd\lambda  \equiv 0$.
The dynamics of the
Reeb vectorfield is very interesting. 
There is a well-known conjecture raised by Arnold in \cite{ar} 
which concerned the Reeb orbit and Legendrian submanifold in 
a contact manifold. If $(\Sigma ,\lambda )$ is a contact manifold 
with contact form $\lambda $ of dimension $2n-1$, then a Legendrian 
submanifold is a submanifold ${\cal L}$ of $\Sigma $, which is 
$(n-1)$dimensional and everywhere tangent to the contact 
structure $\ker \lambda $. Then a characteristic 
chord for $(\lambda ,{{\cal {L}}})$ is a smooth path
$x:[0,T]\to M,T>0$
with
$\dot x(t)=X_{\lambda }(x(t)) \ for \ t\in(0,T),$
$x(0),x(T)\in {\cal {L}}$. 
Arnold raised the following conjectures: 

\vskip 3pt 

{\bf Conjecture1}(see\cite{ar}). Let $\lambda _0$ be the standard tight 
contact form 
$$\lambda _0={{1}\over {2}}(x_1dy_1-y_1dx_1+x_2dy_2-y_2dx_2)$$
on the three sphere 
$$S^3=\{ (x_1,y_1,x_2,y_2)\in R^4|x_1^2+y_1^2+x_2^2+y_2^2=1\}.$$ 
If $f:S^3\to (0,\infty )$ is a smooth function 
and ${\cal {L}}$ is a Legendrian 
knot in $S^3$, then 
there is a characteristic chord for $(f\lambda _0,{\cal {L}})$.

In fact Arnold also conjectured more general cases and multiplicity results just like 
the Lusternik-Schirelman or Morse type number\cite{ar}.

Arnold's conjectures was discussed in \cite{ar,ab, gi}. 
Its solutions 
on the symmetric contact form on $S^3$ and the standard Legendre fibre was given in 
\cite{gi} which also includes multiplicity results.
The complete solution on Conjecture1 was claimed in 1999 in \cite{ma} by using 
the Gromov's nonlinear Fredholm alternative. 
Immediately, the alternate proof was given in \cite{mo}.

\vskip 3pt 

Let $(M,\omega )$ be a symplectic manifold.
Let $J$ be the almost complex structure tamed by $\omega $, i.e., 
$\omega (v, Jv)>0$ for $v\in TM$. Let 
${\cal {J}}$ the space of all tame almost complex structures. 
 
\begin{Definition}
Let 
$$s(M,\omega ,J)=\inf \{ \int _{S^2}f^*\omega >0 |f:S^2\to M \ is \ J-holomorphic \} $$
\end{Definition}

\begin{Definition}
Let 
$$s(M,\omega )=\sup _{J\in {\cal {J}}}l(M,\omega ,J)$$
\end{Definition}

Let $W$ be a Lagrangian submanifold in $M$, i.e., 
$\omega |W=0$. 

\begin{Definition}
Let 
$$l(M,W, \omega )=\inf \{ |\int _{D^2}f^*\omega |>0 |f:(D^2,\partial D^2)\to (M,W) \} $$
\end{Definition}

The main result of this paper is the following:

\begin{Theorem}
Let $(M,\omega )$ be a closed  
compact symplectic manifold or a manifold 
convex at infinity and 
$M\times C$ be a symplectic manifold with 
symplectic form $\omega \oplus \sigma $, 
here $(C,\sigma )$ standard 
symplectic plane. Let 
$2\pi r_0^2<s(M,\omega )$ and 
$B_{r_0}(0)\subset C$ the closed ball with radius 
$r_0$. If 
$(\Sigma ,\lambda )$ be a contact manifold 
of induced type  in $M\times B_{r_0}(0)$ with 
induced contact form $\lambda $, i.e., 
there exists a vector field $X$ transversal to $\Sigma $ such that 
$L_X(\omega \oplus \sigma )=\omega \oplus \sigma $ and $\lambda =i_X(\omega \oplus \sigma )$,  $X_{\lambda } $ its Reeb vector 
field, ${\cal {L}}$ a closed Legendrian submanifold, then 
there exists at least one characteristic chord for 
$(\lambda ,{\cal {L}})$. 
\end{Theorem}
This Theorem generalizes the some results in \cite{ma, mo}. 
For example, if $\omega |\pi _2(M)=0$, then 
$S(M,\omega )=+\infty $. 
We will prove this Theorem by using Lagrangian squeezing theorem which was proved by 
Gromov's nonlinear Fredholm alternative in \cite{ma2} and  
the Mohnke's modification of our Lagrangian construction.

\subsection{Weinstein conjecture}

\begin{Theorem}
Let $(M,\omega )$ be a closed  
compact symplectic manifold or a manifold 
convex at infinity and 
$M\times C$ be a symplectic manifold with 
symplectic form $\omega \oplus \sigma $, 
here $(C,\sigma )$ standard 
symplectic plane. Let 
$2\pi r_0^2<s(M,\omega )$ and 
$B_{r_0}(0)\subset C$ the closed ball with radius 
$r_0$. If 
$(\Sigma ,\lambda )$ be a contact manifold 
of induced type  in $M\times B_{r_0}(0)$ with 
induced contact form $\lambda $, $X_{\lambda } $ its Reeb vector 
field, then 
there exists at least one close characteristics. 
\end{Theorem}
This improves the results in \cite{fhv,hv,ma}.
Again we will prove this Theorem by using Lagrangian squeezing theorem which was proved by 
Gromov's nonlinear Fredholm alternative in \cite{ma2} and  
the Mohnke's modification of our Lagrangian construction.

\section{Lagrangian Squeezing}

\begin{Theorem}
(\cite{ma2})Let $(M,\omega )$ be a closed  
compact symplectic manifold or a manifold convex at infinity and 
$M\times C$ be a symplectic manifold with 
symplectic form $\omega \oplus \sigma $, 
here $(C,\sigma )$ standard 
symplectic plane. Let 
$2\pi r_0^2<s(M,\omega )$ and 
$B_{r_0}(0)\subset C$ the closed disk with radius 
$r_0$. If $W$ is a close Lagrangian manifold in $M\times B_{r_0}(0)$, then 
$$l(M,W, \omega )<2\pi r_0^2$$
\end{Theorem}
This can be considered as an Lagrangian version of Gromov's symplectic squeezing.

\begin{Corollary}
(Gromov\cite{gro})Let $(V',\omega ')$ be an exact symplectic manifold with 
restricted contact boundary and $\omega '=d\alpha ' $. Let 
$V'\times C$ be a symplectic manifold with 
symplectic form $\omega '\oplus \sigma =d\alpha 
=d(\alpha '\oplus \alpha _0$, 
here $(C,\sigma )$ standard 
symplectic plane. 
If $W$ is a close exact Lagrangian submanifold, then 
$l(V'\times C,W, \omega )==\infty $, i.e., there does not exist any 
close exact Lagrangian submanifold in $V'\times C$.
\end{Corollary}

\begin{Corollary}
Let $L^n$ be a close Lagrangian in $R^{2n}$ and $L(R^{2n},L^n,\omega )=2\pi r_0^2>0$, then 
$L^n$ can not be embedded in $B_{r_0}(0)$ as a Lagrangian submanifold. 
\end{Corollary}

\section{Proof Arnold chord conjecture}

\subsection{Constructions of Lagrangian submanifolds}

Let $(\Sigma ,\lambda )$ be a contact manifolds with contact form 
$\lambda $ and $X$ its Reeb vector field, then 
$X$ integrates to a Reeb flow $\eta _t$ for $t\in R^1$. 
Consider the form $d(e^a\lambda )$ 
on the manifold 
$(R\times \Sigma )$, then one can check 
that $d(e^a\lambda )$ is a symplectic 
form on $R\times \Sigma $. Moreover 
One can check that 
\begin{eqnarray}
&&i_X(e^a\lambda )=e^a \\ 
&&i_X(d(e^a\lambda ))=-de^a 
\end{eqnarray}
So, the symplectization of Reeb vector field $X$ is the 
Hamilton vector field of $e^a$ with 
respect to the symplectic form $d(e^a\lambda )$. 
Therefore the Reeb flow lifts to the Hamilton flow 
$h_s$ on $R\times \Sigma $(see\cite{ag,eg}). 

Let ${\cal L}$ be a closed Legendre submanifold 
in $(\Sigma ,\lambda )$, i.e.,  
there exists a smooth embedding $Q:{\cal L}\to \Sigma $ such that 
$Q^*\lambda |_{\cal L}=0$, $\lambda |Q(L)=0$. We also 
write ${\cal {L}}=Q({\cal {L}})$. Let
$$(V',\omega ')=(R\times \Sigma ,d(e^a\lambda ))$$
and 
\begin{eqnarray}
&&W'={\cal {L}}\times R, 
\ \ W'_s={\cal L}\times \{ s\} ;\cr 
&&L'=(0,\cup _s \eta _s(Q({\cal {L}}))), 
\ \ L'_s=(0, \eta _s(Q({\cal {L}})))
\label{eq:3.w}
\end{eqnarray}
define 
\begin{eqnarray}
&&G':W'\to V'  \cr 
&&G'(w')=G'(l,s)=(0,\eta _s(Q(l))) \label{eq:3.ww}
\end{eqnarray}
\begin{Lemma}
There does not exist any Reeb chord connecting Legendre  
submanifold ${\cal {L}}$
in $(\Sigma ,\lambda )$ if and only if  
$G'(W'_s)\cap G'(W'_{s'})$ is empty for $s\ne s'$.
\end{Lemma}
Proof. Obvious. 
\begin{Lemma}
If there does not exist any Reeb chord for $(X_\lambda ,{\cal {L}})$
in $(\Sigma ,\lambda )$ then 
there exists a smooth embedding 
$G':W'\to V'$ with $G'(l,s)=(0,\eta _s(Q(l)))$
such that
\begin{equation}
G'_K:{\cal L}\times (-K, K)\to V'  \label{eq:3.www}
\end{equation}
is a regular open Lagrangian embedding for any finite positive $K$.
We denote $W'(-K,K)=G'_K({\cal L}\times (-K,K))$
\end{Lemma}
Proof. One check 
\begin{equation}
{G'}^*(d(e^a\lambda ))=\eta (\cdot ,\cdot )^*d\lambda 
=(\eta _s^*d\lambda +i_Xd\lambda \wedge ds)=0
\end{equation}
This implies that ${G}'$ is a Lagrangian embedding, this proves 
Lemma3.2. 

\vskip 3pt 

In fact the above proof checks that 
\begin{equation}
{G'}^*(\lambda )=\eta (\cdot , \cdot )^*\lambda =
\eta _s^*\lambda +i_X\lambda ds=ds.
\end{equation}
i.e., $W'$ is an exact Lagrangian submanifold.

\vskip 3pt 

The all above construction was contained in \cite{ma}. Now we intruduce
the Mohnke's upshot. 
Let 
\begin{eqnarray}
&&F':{\cal L}\times R\times R\to R\times \Sigma \cr 
&&F'(l,s,a)=(a,G'(l,s))=(a,\eta _s(Q(l)))  
\end{eqnarray}
Now we embed a elliptic curve $E$ long along $s-axis$ and thin along $a-axis$ such that 
$E\subset [-K,K]\times [0,\varepsilon]$. We parametrize the $E$ by $t$.  
\begin{Lemma}
If there does not exist any Reeb chord for $(X_\lambda ,{\cal {L}})$
in $(\Sigma ,\lambda )$, then 
\begin{eqnarray}
&&F:{\cal L}\times S^1\to R\times \Sigma \cr 
&&F(l,t)=(a(t),G'(l,s(t)))=(a(t),\eta _{s(t)}(Q(l)))  
\end{eqnarray}
is a compact Lagrangian submanifold. Moreover 
\begin{equation}
l(R\times \Sigma ,F({{\cal {L}}}\times S^1,de^a\lambda )=area(E)
\end{equation}
\end{Lemma}
Proof. We check that 
\begin{eqnarray}
F^*(d(e^a\lambda ))
&=&d(F^*(e^{a(t)}\lambda ))\cr
&=&d(e^{a(t))}G'^*\lambda ) \cr
&=&d(e^{a(t)}ds(t)) \cr
&=&e^{a(t)}(a_tdt\wedge s_tdt)\cr
&=&0
\end{eqnarray}
which shows that $F$ is a Lagrangian embedding. 

If the circle $C$ homotopic to $C_1\subset {\cal {L}}\times s_0$ then  we compute
\begin{eqnarray}
\int _CF^*(e^a\lambda )=\int _{C_1}F^*(e^a\lambda )=0. 
\end{eqnarray}
since $\lambda |C_1=0$ due to $C_1\subset {\cal {L}}$ and 
$\cal L$ is Legendre submanifold. 

If the circle $C$ homotopic to $C_1\subset l_0\times S^1$ then  we compute
\begin{eqnarray}
\int _CF^*(e^a\lambda )=\int _{C_1}F^*(e^a\lambda )=n(area(E)). 
\end{eqnarray}
This proves the Lemma.

\subsection{Proof on Theorem 1.1}

Since $(\Sigma ,\lambda )$ be a contact manifold 
of induced type  in $M\times B_{r_0}(0)$ with 
induced contact form $\lambda $, then by the well known theorem that 
the neighbourhood $(U(\Sigma ),\omega )$ of $\Sigma $ is symplectomorphic to 
$([-\varepsilon ,\varepsilon]\times \Sigma ,de^a\lambda )$ for small 
$\varepsilon $. So, by Lemma 3.3, we have 
a close Lagrangian submanifold $F({\cal {L}}\times S^1)$ contained in 
$M\times B_{r_0}(0)$. By Lagrangian squeezing theorem, i.e., Theorem 2.1, 
we have 
\begin{equation}
l(M\times C,F({{\cal {L}}}\times S^1,\omega )=area(E)\leq 2\pi r_0^2.
\end{equation}
If $K$ large enough, $area(E)>2\pi r_0^2$. This 
is a contradiction. This contradiction shows 
there exists at least one characteristic chord for 
$(\lambda ,{\cal {L}})$.

\section{Proof on Weinstein conjecture}

\subsection{Constructions of Lagrangian submanifolds}

Let $(\Sigma ,\lambda )$ be a contact manifolds with contact form 
$\lambda $ and $X$ its Reeb vector field, then 
$X$ integrates to a Reeb flow $\eta _t$ for $t\in R^1$. 
Let 
$$(V',\omega ')=((R\times \Sigma )\times (R\times \Sigma ), 
d(e^a\lambda )\ominus d(e^b\lambda ))$$ 
and 
$${\cal {L}}=\{ ((0,\sigma ),(0,\sigma ))|(0,\sigma )\in R\times \Sigma \}.$$
Let 
$$L'={\cal {L}}\times R,
L_s'={\cal {L}}\times \{s\}.$$
Then 
define 
\begin{eqnarray}
&&G':L'\to V'\cr 
&&G'(l')=G'(((\sigma ,0),(\sigma ,0)),s)
=((0,\sigma ),(0,\eta _s(\sigma )))
\end{eqnarray}
Then
    
$$W'=G'(L')=\{ ((0,\sigma ),(0,\eta _s(\sigma )))
|(0,\sigma )\in R\times \Sigma ,s\in R\}$$

$$W_s'=G'(L'_s)=\{ ((0,\sigma ),(0,\eta _s(\sigma )))
|(0,\sigma )\in R\times \Sigma \}$$
for fixed $s\in R$.

\begin{Lemma}
There does not exist any Reeb closed orbit in 
$(\Sigma ,\lambda )$ if and only if  
$W'_s\cap W'_{s'}$ is empty for $s\ne s'$.
\end{Lemma}
Proof. First if there exists a closed Reeb orbit in 
$(\Sigma ,\lambda )$, i.e., there exists 
$\sigma _0\in \Sigma $, $t_0>0$ such that 
$\sigma _0=\eta _{t_0}(\sigma _0)$, then 
$((0,\sigma _0),(0,\sigma _0))\in W'_0\cap W'_{t_0}$.  
Second if there exists $s_0\ne s_0'$ such 
that $W'_{s_0}\cap W'_{s_0'}\ne \emptyset $, i.e., 
there exists $\sigma _0$ such that 
$$((0,\sigma _0),(0,\eta _{s_0}(\sigma _0))
=((0,\sigma _0),(0,\eta _{s_0'}(\sigma _0)),$$ 
then $\eta _{(s_0-s_0')}(\sigma _0)=\sigma _0$, i.e., 
$\eta _t(\sigma _0)$ is a closed Reeb orbit.

\begin{Lemma}
If there does not exist any closed Reeb orbit in 
$(\Sigma ,\lambda )$ then 
there exists a smooth Lagrangian injective immersion 
$G':W'\to V'$ with $G'(((0,\sigma ),(0,\sigma )),s)
=((0,\sigma ),(0,\eta _s(\sigma )))$
such that
\begin{equation}
G'_{s_1,s_2}:{\cal L}\times (-s_1,s_2)\to V'
\end{equation}
is a regular exact Lagrangian embedding for any finite real number 
$s_1$, $s_2$, here we denote by $W'(s_1,s_2)=G'_{s_1,s_2}({\cal L}\times 
(s_1,s_2))$. 
\end{Lemma}
Proof. One check 
\begin{equation}
{G'}^*((e^a\lambda -e^b\lambda ))
=\lambda -\eta (\cdot ,\cdot )^*\lambda 
=\lambda -(\eta _s^*\lambda +i_X\lambda ds)=-ds
\end{equation}
since $\eta _s^*\lambda =\lambda $. 
This implies that ${G}'$ is an exact  Lagrangian embedding, this proves 
Lemma 3.2. 

\vskip 3pt 

Now we modify the above construction as follows:  
\begin{eqnarray}
&&F':{\cal {L}}\times R\times R\to (R\times \Sigma )\times 
(R\times \Sigma )\cr 
&&F'(((0,\sigma ),(0,\sigma )),s,b)=((0,\sigma ),(b,\eta _s(\sigma )))
\end{eqnarray}
Now we embed a elliptic curve $E$ long along $s-axis$ and thin along $b-axis$ such that 
$E\subset [-s_1,s_2]\times [0,\varepsilon]$. We parametrize the $E$ by $t$.

\begin{Lemma}
If there does not exist any closed Reeb orbit 
in $(\Sigma ,\lambda )$, 
then 
\begin{eqnarray}
&&F:{\cal {L}}\times S^1\to (R\times \Sigma )\times 
(R\times \Sigma )\cr 
&&F(((0,\sigma ),(0,\sigma )),t)=((0,\sigma ),(b(t),\eta _{s(t)}(\sigma )))
\end{eqnarray}
is a compact Lagrangian submanifold. Moreover 
\begin{equation}
l(V',F({{\cal {L}}}\times S^1,d(e^a\lambda -e^b\lambda ))=area(E)
\end{equation}
\end{Lemma}
Proof. We check
that 
\begin{eqnarray}
{F}^*(e^a\lambda \ominus e^b\lambda )&=&-e^{b(t)}ds(t)
\end{eqnarray}
So, $F$ is a Lagrangian embedding.

If the circle $C$ homotopic to $C_1\subset {\cal {L}}\times s_0$ then  we compute
\begin{eqnarray}
\int _CF^*(e^a\lambda )=\int _{C_1}F^*(e^a\lambda )=0. 
\end{eqnarray}
since $\lambda |C_1=0$ due to $C_1\subset {\cal {L}}$ and 
$\cal L$ is Legendre submanifold. 

If the circle $C$ homotopic to $C_1\subset l_0\times S^1$ then  we compute
\begin{eqnarray}
\int _CF^*(e^a\lambda )=\int _{C_1}F^*(e^a\lambda )=n(area(E)). 
\end{eqnarray}
This proves the Lemma.

\subsection{Proof on Theorem 1.2}

Since $(\Sigma ,\lambda )$ be a contact manifold 
of induced type  in $M\times B_{r_0}(0)$ with 
induced contact form $\lambda $, then by the well known theorem that 
the neighbourhood $(U(\Sigma ),\omega )$ of $\Sigma $ is symplectomorphic to 
$([-\varepsilon ,\varepsilon]\times \Sigma ,de^a\lambda )$ for small 
$\varepsilon $. So, by Lemma 4.3, we have 
a close Lagrangian submanifold $F({\cal {L}}\times S^1)$ contained in 
$M\times C\times M\times B_{r_0}(0)$. By Lagrangian squeezing theorem, i.e., Theorem 2.1, 
we have 
\begin{equation}
l((M\times C)\times (M\times C),F({{\cal {L}}}\times S^1,\omega \oplus \omega )=area(E)\leq 2\pi r_0^2.
\end{equation}
If $s_2-s_1$ large enough, $area(E)>2\pi r_0^2$. This 
is a contradiction. This contradiction shows 
there exists at least one close characteristics.

\end{document}